\documentclass[number,dvips]{arxbj}
\usepackage{cursive}  

\aid{0}
\volume{16}
\issue{2}
\pubyear{2010}
\firstpage{331}
\lastpage{342}
\doi{10.3150/09-BEJ214}

\makeatletter

\newcommand{\R}{\mathbb{R}}
\newcommand{\Prob}{\mathbb{P}}
\newcommand{\Z}{\mathbb{Z}}
\newcommand{\Ex}{\mathbb{E}}
\newcommand{\eqref}[1]{(\ref{#1})}

\newproclaim{definition}{Definition}

\newremark{remark}{Remark}
\newremark{example}{Example}
\makeatother

\begin{document}
\begin{frontmatter}

\title{Copulas for Markovian dependence}
\runtitle{Copulas for Markovian dependence}

\begin{aug}
\author{\fnms{Andreas N.} \snm{Lager{\aa}s}\ead[label=e1]{andreas@math.su.se}\corref{}}
\runauthor{A.N. Lager\aa s}
\pdfauthor{Andreas N. Lageras}
\address{Department of Mathematics, Stockholm University, SE-10691
Stockholm, Sweden.\\ \printead{e1}}
\end{aug}

\received{\smonth{12} \syear{2008}}
\revised{\smonth{4} \syear{2009}}

%
\begin{abstract}
Copulas have been popular to model dependence for multivariate
distributions, but have not been used much in modelling temporal
dependence of univariate time series. This paper demonstrates some
difficulties with using copulas even for Markov processes: some
tractable copulas such as mixtures between copulas of complete co- and
countermonotonicity and independence (Fr{\'{e}}chet copulas) are shown
to imply quite a restricted type of Markov process and Archimedean
copulas are shown to be incompatible with Markov chains. We also
investigate Markov chains that are spreadable or, equivalently,
conditionally i.i.d.
\end{abstract}

%
\begin{keyword}
\kwd{copulas}
\kwd{exchangeability}
\kwd{Markov chain}
\kwd{Markov process}
\end{keyword}

\end{frontmatter}

\section{Introduction}

Copulas, which will be defined in Section \ref{intro_sec}, describe the
dependence of a multivariate distribution that is invariant under
monotone (increasing) transformations of each coordinate. In this
paper, we investigate the dependence that arises in a one-dimensional
Markov process. Darsow \textit{et al.} \cite{dno} began the study of
copulas related to Markov processes; see also \cite{nelsen}, Chapter
6.3. More precisely, they showed what the Kolmogorov--Chapman equations
for transition kernels translate to in the language of copulas and
introduced some families of copulas $(C_{st})_{s\leq t}$ that are
consistent in the sense that $C_{st}$ is the copula of $(X_s,X_t)$ for
a Markov process $(X_t)_{t\geq0}$.

In Section \ref{intro_sec}, we will introduce a Markov product of
copulas $C\ast D$ such that if $C$ gives the dependence of $(X_0,X_1)$
and $D$ the dependence of $(X_1,X_2)$, then $C\ast D$ gives the
dependence of $(X_0,X_2)$ for a Markov chain $X_0,X_1,X_2$. An analogy
is that of a product of transition matrices of finite-state Markov
chains, in particular, doubly stochastic matrices (whose column sums
are all 1) since they have uniform stationary distribution.

This approach might, at first, seem like a sensible way of introducing
the machinery of copulas into the field of stochastic processes:
Mikosch \cite{mikosch}, for example, has criticized the widespread use
of copulas in many areas and, among other things, pointed out a lack of
understanding of the temporal dependence, in terms of copulas, of most
basic stochastic processes.

This paper builds on that of Darsow \textit{et al.}~\cite{dno}, but
with a heavier emphasis on probabilistic, rather than analytic or
algebraic, understanding. Our main results are negative, in that we
show how:
\begin{enumerate}[(3)]
\item[(1)] a proposed characterization of the copulas of time-homogeneous
Markov processes fails (Section \ref{families_sec});
\item[(2)] Fr{\'{e}}chet copulas imply quite strange Markov processes
(Section \ref{frechet_sec});
\item[(3)] Archimedean copulas are incompatible with the dependence of
Markov chains (Section~\ref{arch_sec});
\item[(4)] a conjectured characterization of idempotent copulas,
related to
exchangeable Markov chains, fails (Section \ref{idem_sec}).
\end{enumerate}

\section{Copulas and the Markov product}\label{intro_sec}

\begin{definition}
A \textit{copula} is a distribution function of a multivariate random
variable whose univariate marginal distributions are all uniform on $[0,1]$.
\end{definition}

We will mostly concern ourselves with two-dimensional copulas. In the
following, all random variables denoted by $U$ have a uniform
distribution on $[0,1]$ (or, sometimes, $(0,1)$).

\begin{definition}
$\Pi(x,y)=xy$ is the copula of independence: $(U_1,U_2)$ has the
distribution $\Pi$ if and only if $U_1$ and $U_2$ are independent.
\end{definition}

\begin{definition}
$M(x,y)=\min(x,y)$ is the copula of complete positive dependence:
$(U_1,U_2)$ has the distribution $M$ if and only if $U_1=U_2$ almost
surely (a.s.).
\end{definition}

\begin{definition}
$W(x,y)=\max(x+y-1,0)$ is the copula of complete negative dependence:
$(U_1,U_2)$ has the distribution $W$ if and only if $U_1=1-U_2$ a.s.
\end{definition}

Note that a mixture $\sum_ip_i C_i$ of copulas $C_1,C_2,\ldots$ is also
a copula if $p_1,p_2,\ldots$ is a probability distribution since one can
interpret the mixture as a randomization: first, choose a copula
according to the distribution $p_1,p_2,\dots$ and then draw from the
chosen distribution.

It is well known that if $X$ is a (one-dimensional) continuous random
variable with distribution function $F$, then $F(X)$ is uniform on
$[0,1]$. Thus, if $(X_1,\dots,X_n)$ is an $n$-dimensional continuous
random variable with joint distribution function $F$ and marginal
distributions $(F_1,\dots,F_n)$, then the random variable
$(F_1(X_1),\dots,F_n(X_n))$ has uniform marginal distributions, that
is,\ its joint distribution function is a copula, say $C$.

\textit{Sklar's theorem} (see \cite{nelsen}, Theorem 2.10.9) states that
any $n$-dimensional distribution function~$F$ with marginals
$(F_1,\dots
,F_n)$ can be ``factored'' into $F(x_1,\dots,x_n)=C(F_1(x_1),\dots
, F_n(x_n))$ for a copula $C$, which is, furthermore, unique if the
distribution $F$ is continuous. We say that the $n$-dimensional
distribution $F$, or the random variable $(X_1,\dots,X_n)$, has the
copula $C$.

\begin{remark}
When $(X_1,\dots,X_n)$ does not have a unique copula, all copulas of
this random variable agree at points $(u_1,\dots,u_n)$ where $u_i$ is
in the range $R_i$ of the function $x_i\mapsto F_i(x_i)$. One can
obtain a unique copula by an interpolation between these points which
is linear in each coordinate, and we will, as Darsow \textit{et al.}\
\cite{dno}, speak of this as the copula of such random variables.
\end{remark}

Copulas allow for a study of the dependence in a multivariate
distribution separately from the marginal distributions. It gives
reasonable information about dependence in the sense that the copula is
unchanged if $(X_1,\dots,X_n)$ is transformed into $(g_1(X_1),\dots
,g_n(X_n))$, where $g_1,\dots,g_n$ are strictly increasing.

\begin{example}
The notion of copulas makes it possible to take a copula from, say, a
multivariate $t$-distribution and marginal distributions from, say, a
normal distribution and combine them into a multivariate distribution
where the marginals are normal, but the joint distribution is not
multivariate normal. This is sometimes desirable in order to have
models with, in a sense, ``stronger'' dependence than what is possible
for a multivariate normal distribution.
\end{example}

\begin{example}
$(X_1,X_2)$ has the copula $\Pi$ if and only if $X_1$ and $X_2$ are
\textit{independent}. $(X_1,X_2)$ has the copula $M$ if and only if
$X_2=g(X_1)$ for a strictly \textit{increasing} function $g$. $(X_1,X_2)$
has the copula $W$ if and only if $X_2=h(X_1)$ for a strictly \emph
{decreasing} function $h$. (When $X_1$ and $X_2$ furthermore have the
same marginal distributions, they are usually called \textit{antithetic}
random variables.)
\end{example}

In this paper, we are in particular interested in the dependence that
arises in a Markov process in $\R$, for example,\ the copula of
$(X_0,X_1)$ for a stationary Markov chain $X_0,X_1,\ldots.$ By \cite
{kallenberg}, Proposition 8.6, the sequence $X_0,X_1,\ldots$ constitutes
a Markov chain in a Borel space $S$ if and only if there exist
measurable functions $f_1,f_2,\ldots\dvtx S\times[0,1]\to S$ and i.i.d.
random variables $V_1,V_2,\dots$ uniform on $[0,1]$ and all independent
of $X_0$ such that $X_n=f_n(X_{n-1},V_n)$ a.s. for $n=1,2,\ldots.$ One
may let $f_1=f_2=\cdots= f$ if and only if the process is time-homogeneous.

We can, without loss of generality, let $S=[0,1]$ since we can
transform the coordinates $X_0,X_1,\dots$ monotonically without
changing their copula. The copula is clearly related to the function
$f$ above. We have $f_{\Pi}(x,u)=u$, $f_M(x,u)=x$ and $f_W(x,u)=1-x$
with obvious notation.

Darsow \textit{et al.}\ \cite{dno} introduced an operation on copulas
denoted $\ast$ which we will call the \textit{Markov product}.

\begin{definition}
Let $X_0,X_1,X_2$ be a Markov chain and let $C$ be the copula of
$(X_0,X_1)$, $D$ the copula of $(X_1,X_2)$ and $E$ the copula of
$(X_0,X_2)$ (note that $X_0,X_2$ is also a Markov chain). We then write
$C\ast D=E$.
\end{definition}

It is also possible to define this operation as an integral of a
product of partial derivatives of the copulas $C$ and $D$;
see \cite{dno}, formula (2.10), or \cite{nelsen}, formula (6.3.2),
but, in this paper, the
probabilistic definition will suffice.

From the definition, it should be clear that the operation $\ast$ is
associative, but not necessarily commutative and for all $C$:
\begin{eqnarray*}
\Pi\ast C &=& C\ast\Pi= \Pi,\\
M\ast C &=& C\ast M = C
\end{eqnarray*}
so that $\Pi$ acts as a null element and $M$ as an identity. We write
$C^{\ast n}$ for the $n$-fold Markov product of $C$ with itself and
define $C^{\ast0}=M$. We have $W^{\ast2}=M$, so $W^{\ast n}=M$ if $n$
is even and $W^{\ast n}=W$ is $n$ is odd. In Section \ref{idem_sec}, we
will investigate \textit{idempotent} copulas $C$, meaning $C^{\ast2}=C$.

\begin{example}
If $X_0,X_1,\dots$ is a time-homogeneous Markov chain where $(X_0,X_1)$
has copula $C$, then $C^{\ast n}$ is the copula of $(X_0,X_n)$ for all
$n=0,1,\dots.$
\end{example}

\begin{definition}
For any copula $C(x,y)$ of the random variable $(X,Y)$, we define its
\textit{transpose} $C^T(x,y)=C(y,x)$, the copula of $(Y,X)$.
\end{definition}

We can say that $W$ is its own \textit{inverse} since $W\ast W=M$.
\begin{definition}
In general, we say that a copula $R$ is \textit{left-invertible} or a
\textit{right-inverse} if there exists a copula $L$ such that $L\ast R=M$
and we say that $L$ is \textit{right-invertible} or a \textit{left-inverse}.
\end{definition}

The equation $L\ast R=M$ implies that any randomness in the transition
described by $L$ is eliminated by $R$ and thus $f_R(x,u)$ must be a
function of $x$ alone. A rigorous proof of the last proposition may be
found in \cite{dno}, Theorem 11.1. Furthermore, if $L$ is a
right-invertible copula of $(X,Y)$, then its right-inverse $R$ can be
taken as the transpose of $L$, $R=L^{\mathrm{T}}$, since $M$ is the copula of
$(X,X)$ and thus $R$ should be the copula of $(Y,X)$ so that $L,R$
correspond to the Markov chain $X,Y,X$. A proof of this can also found
in \cite{dno}, Theorem 7.1.

\begin{example}
Let $L_{\theta}$ be the copula of the random variable $(X,Y)$ whose
distribution is as follows: $(X,Y)$ is uniform on the line segment
$y=\theta x$, $0\leq x\leq1$, with probability $0\leq\theta\leq1$
and $(X,Y)$ is uniform on the line segment $y=1-(1-\theta)x$, $0\leq
x\leq1$, with probability $1-\theta$. The function
\[
f_{L_{\theta}}(x,u)=\theta x\mathbf{1}(u\leq\theta)+\bigl
(1-(1-\theta
)x\bigr)\mathbf{1}(u>\theta)
\]
can be used to describe the transition from $X$ to $Y$. Let $R_{\theta
}=L_{\theta}^{\mathrm{T}}$. One can take
\[
f_{R_{\theta}}(y,v)=\frac{y}{\theta}\mathbf{1}(y\leq\theta)+\frac
{1-y}{1-\theta}\mathbf{1}(y>\theta)
\]
to describe the transition from $Y$ to $X$. Note that $f_{R_{\theta
}}(y,v)$ is a function of $y$ only. We also get $f_{R_{\theta
}}(f_{L_{\theta}}(x,u),v)=f_M(x,w)=x$ so that, indeed, $L_{\theta
}\ast
R_{\theta} = M$.
\end{example}

The Markov product is linear:
\begin{equation}\label{linear}
\sum_ip_iC_i \ast\sum_jq_j D_j = \sum_{ij}p_iq_j C_i\ast D_j
\end{equation}
since the left-hand side can be interpreted as first choosing a $C_i$
with probability $p_i$ as transition mechanism from $X_0$ to $X_1$ and
then independently choosing a $D_j$ with probability $q_j$ as
transition mechanism from $X_1$ to $X_2$, whereas the right-hand side
can be interpreted as choosing a combined transition mechanism $C_i\ast
D_j$ from $X_0$ to $X_2$ with probability $p_iq_j$.

For a given Markov process $(X_t)_{t\geq0}$ in continuous time, we
will denote the copula of $(X_s,X_t)$ by $C_{st}$ for $s\leq t$. For
time-homogeneous processes, we only write $C_t$ for the copula of
$(X_s,X_{s+t})$ for $t\geq0$. Note that, for all $t$,
\[
C_{tt} = C_{00} = C_0 = M.
\]
Copulas for continuous-time Markov processes must obey a
Kolmogorov--Chapman-type relationship:
\begin{equation}\label{kol-chap}
C_{rt}=C_{rs}\ast C_{st},\qquad r\leq s\leq t.
\end{equation}

\section{Some families of copulas}\label{families_sec}

Let $(X_t)_{t\geq0}$ be a Markov process with transition kernel
$P_{st}(x,\cdot)$ and marginal distributions $(F_t)_{t\geq0}$. Now,
\begin{equation}\label{krangligt}
C_{st}(F_s(x),F_t(y))=\Prob(X_s\leq x,X_t\leq y)=\int_{-\infty
}^xP_{st}(u,(-\infty,y])\,\mathrm{d}F_s(u)
\end{equation}
and, from this, $C_{st}$ may be derived in principle.

The expression \eqref{krangligt} becomes more manageable if the
marginal distributions are uniform and if the transition kernel
furthermore has a density $f_{st}(x,y)$, then we get that the density
$c_{st}(x,y)=\frac{\partial^2}{\partial x\,\partial y}C_{st}(x,y)$ of the
Markov copula equals the transition density: $c_{st}=f_{st}$.

\begin{example}\label{ex_tva}
Let $(U_t)_{t\geq0}$ be a Brownian motion reflected at 0 and 1, with
$\sigma=1$ and with $U_0$ uniform on $(0,1)$. This process is
stationary and time-homogeneous, with
\[
c_t(x,y)=\frac{1}{\sqrt{2\curpi t}}\sum_{n\in\Z}\bigl(\mathrm{e}^{-
(2n+y-x)^2/(2t)}+\mathrm{e}^{-(2n-y-x)^2/(2t)}\bigr).
\]
It is clear that $C_t\to M$ as $t\to0$ and $C_t\to\Pi$ as $t\to
\infty$.
\end{example}

It is usually hard to compute transition densities for interesting
processes, so another way of obtaining families of Markov copulas is to
construct them directly from copulas so that \eqref{kol-chap} holds. A
problem with this approach is that a probabilistic understanding of the
process may be lost.

\begin{example}
Darsow \textit{et al.} \cite{dno} pose the question of
whether all time-homo\-geneous Markov copulas may be expressed as
\begin{equation}\label{inte_all_hom}
C_t = \mathrm{e}^{-at}\Biggl(E+\sum_{n=1}^{\infty}\frac{a^nt^n}{n!}C^{\ast
n}\Biggr),
\end{equation}
where $a$ is a positive constant and $E$ and $C$ are two copulas
satisfying $C\ast E=E\ast C=C$ and $E$ is idempotent ($E\ast E=E$). We
immediately observe that $C_0=E$ according to equation \eqref
{inte_all_hom} and thus $E$ cannot be taken to be arbitrary, but must
equal $M$. However, since $M$ commutes with all copulas, $C$ may be
arbitrary. As $M=C^{\ast0}$, we can rewrite
\begin{equation}\label{inte_all_hom_po}
C_t = \sum_{n=0}^{\infty}\frac{(at)^n}{n!}\mathrm{e}^{-at}C^{\ast n} = \Ex
\bigl[C^{\ast N(t)}\bigr],
\end{equation}
where $N$ is a Poisson process with intensity $a$. We can thus give the
following probabilistic interpretation: a Markov process has the Markov
copula of equation \eqref{inte_all_hom_po} if it jumps according to the
Poisson process $N$ with intensity $a$ and, at each jump, it jumps
according to the copula $C$. Between jumps, it remains constant. This
clearly does not cover all possible time-homogeneous Markov processes
or Markov copulas; see\ the previous Example \ref{ex_tva}.
\end{example}

\section{Fr{\'{e}}chet copulas}\label{frechet_sec}

In this section, we only consider Markov processes in continuous time.

A copula $C$ is said to be in the Fr{\'{e}}chet family if $C=\alpha W +
(1-\alpha-\beta)\Pi+ \beta M$ for some non-negative constants
$\alpha$
and $\beta$ satisfying $\alpha+\beta\leq1$; see \cite{nelsen}, page 12.
Darsow \textit{et al.} \cite{dno} found conditions on the functions
$\alpha(s,t)$ and $\beta(s,t)$ in
\[
C_{st}=\alpha(s,t) W + \bigl(1-\alpha(s,t)-\beta(s,t)\bigr)\Pi+ \beta(s,t) M
\]
such that $C_{st}$ satisfies equation \eqref{kol-chap}. By equation
\eqref{linear}, we find
\begin{eqnarray}
\label{alpha_eq}
\beta(r,s)\alpha(s,t)+\alpha(r,s)\beta(s,t)&=&\alpha(r,t),\\
\label{beta_eq}
\alpha(r,s)\alpha(s,t)+\beta(r,s)\beta(s,t)&=&\beta(r,t).
\end{eqnarray}
Darsow \textit{et al.} \cite{dno} solved these equations by putting
$r=0$ and defining $f(t)=\alpha(0,t)$ and $g(t)=\beta(0,t),$ which yields
\begin{eqnarray*}
\alpha(s,t)&=&\frac{f(t)g(s)-f(s)g(t)}{g(s)^2-f(s)^2},\\
\beta(s,t) &=&\frac{g(t)g(s)-f(s)f(t)}{g(s)^2-f(s)^2}.
\end{eqnarray*}
This solution in terms of the functions $f$ and $g$ does not have an
immediate probabilistic interpretation and it is therefore hard to give
necessary conditions on the functions $f$ and $g$ for \eqref{alpha_eq}
and \eqref{beta_eq} to hold.

We will first investigate the time-homogeneous case, where $\alpha
(s,t)=a(t-s)$ and $\beta(s,t)=b(t-s)$ for some functions $a$ and $b$.
The equations \eqref{alpha_eq} and \eqref{beta_eq} are then
\begin{eqnarray}
\label{alpha_hom}
b(s)a(t)+a(s)b(t)&=&a(s+t),\\
\label{beta_hom}
a(s)a(t)+b(s)b(t)&=&b(s+t).
\end{eqnarray}
Letting $\rho(t)=a(t)+b(t)$, we find, by summing the two equations
\eqref{alpha_hom} and \eqref{beta_hom}, that
\begin{equation}\label{rho_hom}
\rho(s)\rho(t)=\rho(s+t).
\end{equation}
Since $\rho$ is bounded and $\rho(0)=1$ (since $C_0=M$), we necessarily
have $\rho(t)=\mathrm{e}^{-\lambda t}$, where $\lambda\geq0$ or $\rho
(t)=\mathbf
{1}(t=0)$. Note that $\rho(t)$ equals the probability that a Poisson
process $N_{\Pi}$ with intensity~$\lambda$ has no points in the
interval $(0,t]$.

For the moment, we disregard the possibility $\rho(t)=\mathbf{1}(t=0)$.
Since $\rho$ is positive, we can define $\sigma(t)=a(t)/\rho(t)$. By
dividing both sides of \eqref{beta_hom} by $\rho(s+t)$ and using
\eqref{rho_hom}, we get
\begin{equation}\label{sigma_hom}
\sigma(s)\sigma(t)+\bigl(1-\sigma(s)\bigr)\bigl(1-\sigma(t)\bigr)=\sigma(s+t).
\end{equation}
If we now let $\tau(t)=1-2\sigma(t)$, equation \eqref{sigma_hom} yields
\begin{equation}
\tau(s)\tau(t)=\tau(s+t)
\end{equation}
and, by the same reasoning as for $\rho$, we get $\tau(t)=\mathrm{e}^{-2\mu t}$
for some $\mu\geq0$ or $\tau(t)=\mathbf{1}(t=0)$. We disregard the
latter possibility for the moment. Thus, $\sigma(t)=\frac12-\frac
12\mathrm{e}^{-2\mu t}=\mathrm{e}^{-\mu t}\sinh\mu t$ for some constant $\mu\geq0$.
Note that
\begin{equation}\label{odd}
\sigma(t)=\mathrm{e}^{-\mu t}\sinh\mu t = \sum_{k=0}\mathrm{e}^{-\mu t}\frac{(\mu
t)^{2k+1}}{(2k+1)!},
\end{equation}
that is,\ $\sigma(t)$ equals the probability that a Poisson process
$N_{W}$ with intensity $\mu$ has an odd number of points in $(0,t]$.

Thus, we have
\begin{eqnarray}\label{hom}
C_t &=& \sigma(t)\rho(t)W+\bigl(1-\rho(t)\bigr)\Pi+\bigl(1-\sigma(t)\bigr)\rho(t)M
\nonumber\\
&=& \mathrm{e}^{-(\lambda+\mu)t}\sinh(\mu t) W+(1-\mathrm{e}^{-\lambda t})\Pi+\mathrm{e}^{-(\lambda
+\mu)t}\cosh(\mu t) M\nonumber\\[-8pt]\\[-8pt]
&=&\Prob\bigl(N_W(t)\mbox{ is odd}, N_{\Pi}(t)=0\bigr)W+\Prob\bigl(N_{\Pi}(t)\geq1\bigr)\Pi\nonumber\\
&&{}+\Prob\bigl(N_W(t)\mbox{ is even}, N_{\Pi}(t)=0\bigr)M,\nonumber
\end{eqnarray}
where the aforementioned Poisson processes $N_{\Pi}$ and $N_W$ are independent.

\subsection*{Probabilistic interpretation}

The time-homogeneous Markov
process with $C_t$ as copula is therefore rather special. We may,
without loss of generality, assume that all marginal distributions are
uniform on $[0,1]$. It ``restarts'' -- becoming independent of its
history -- according to a Poisson process $N_{\Pi}$. Independently of
this process, it ``switches'' by transforming a present value $U_{t-}$
to $U_t=1-U_{t-}$, and this happens according to a Poisson process
$N_W$. Note that the intensity of either process may be zero.

If $\tau(t)=\mathbf{1}(t=0)$, then $\sigma(t)=\frac12\mathbf{1}(t>0)$
so that $C_t = \rho(t)(\frac12W+\frac12M)+(1-\rho(t))\Pi$ for \mbox{$t>0$}.
The process can be described as follows. Between points $t_i<t_{i+1}$
of $N_{\Pi}$, the random variables $(U_t)_{t_i\leq t<t_{i+1}}$ are
independent and have the distribution
$P(U_t=U_{t_i})=P(U_t=1-U_{t_i})=\frac12$.

If $\rho(t)=\mathbf{1}(t=0)$, then we have $C_t = \Pi$ for $t>0$ so
that the process is, at each moment, independent of the value at any
other moment, that is,\ $(U_t)_{t\geq0}$ is a collection of
independent random variables.

With the probabilistic interpretation, it is easy to rewrite equation
\eqref{hom} in the form \eqref{inte_all_hom_po} when $\rho,\sigma>0$.
The process makes a jump of either ``restart'' or ``switch'' type with
intensity $\lambda+\mu$ and each jump is of ``restart'' type with
probability $\lambda/(\lambda+\mu)$ and of ``switch'' type with
probability $\mu/(\lambda+\mu)$. Thus,
\[
C_t = \sum_{n=0}^{\infty}\frac{((\lambda+\mu)t)^n}{n!}\mathrm{e}^{-(\lambda
+\mu
)t}\biggl(\frac{\lambda}{\lambda+\mu}\Pi+\frac{\mu}{\lambda+\mu
}W\biggr)^{\ast n} = \Ex\bigl[C^{\ast N(t)}\bigr],
\]
where $C=\frac{\lambda}{\lambda+\mu}\Pi+\frac{\mu}{\lambda+\mu
}W$ and
$N$ is a Poisson process with intensity $\lambda+\mu$.

It is clear that the time-homogeneous process can be generalized to a
time-inhomogeneous Markov process by taking $N_{\Pi}$ and $N_W$ to be
independent inhomogeneous Poisson processes. With
\begin{eqnarray*}
\rho(s,t)&=&\Prob\bigl(N_{\Pi}(t)-N_{\Pi}(s)=0\bigr),\\
\sigma(s,t)&=&\Prob\bigl(N_W(t)-N_W(s)\mbox{ is odd}\bigr),
\end{eqnarray*}
we get a more general version of the Fr{\'{e}}chet copula:
\[
C_{st}=\sigma(s,t)\rho(s,t)W+\bigl(1-\rho(s,t)\bigr)\Pi+\bigl(1-\sigma(s,t)\bigr)\rho(s,t)M,
\]
with essentially the same probabilistic interpretation as the
time-homogeneous case.

In the time-inhomogeneous case, it is also possible to let either or
both of the two processes consist of only one point, say $\tau_{\Pi}$
and/or $\tau_{W}$ that may have arbitrary distributions on $(0,\infty
)$. In addition to this, both in the Poisson case and the single point
case, it is possible to add deterministic points to the processes
$N_{\Pi}$ and $N_{W}$ and still retain the (time-inhomogeneous) Markov
property.

\section{Archimedean copulas}\label{arch_sec}

If a copula of an $n$-dimensional random variable $(X_1,\dots,X_n)$ is
of the form
\begin{equation}\label{arch}
\phi^{-1}\bigl(\phi(u_1)+\cdots+\phi(u_n)\bigr),
\end{equation}
it is called \textit{Archimedean} with generator $\phi$. 
Necessary and sufficient conditions on the generator to produce a
copula are given in \cite{mn}. We will only use the following necessary
properties, which we express with the inverse $\psi=\phi^{-1}$. The
function $\psi$ is non-increasing, continuous, defined on $[0,\infty)$
with $\psi(0)=1$ and $\lim_{x\to\infty}\psi(x)=0$, and decreasing when
$\psi>0$; see \cite{mn}, Definition~2.2.

\begin{example}\label{arch_ex}
Let $(U_1,\dots,U_n)$ be distributed according to \eqref{arch} and let
$X_i = -\phi(U_i)$ for $i=1,\dots,n$. Since $-\phi$ is an increasing
function, $(X_1,\dots,X_n)$ has the same copula as $(U_1,\dots,U_n)$.
All components of $(X_1,\dots,X_n)$ are non-positive and
\begin{eqnarray*}
\Prob(X_i\leq-x_i, i=1,\dots,n)&=&\Prob\bigl(-\phi(U_i)\leq
-x_i,i=1,\dots
,n\bigr)\\
&=&\Prob\bigl(U_i\leq\psi(x_i), i=1,\dots,n\bigr) \\
&=& \psi(x_1+\cdots+x_n)
\end{eqnarray*}
for $x_1,\dots,x_n\geq0$.
\end{example}

Let us now consider Markov processes with Archimedean copulas.
\begin{proposition*}
If $X_1,\dots,X_n$ is a Markov chain, where $(X_1,\dots,X_n)$ has an
Archimedean copula with generator $\phi$ and $n\geq3$, then all
$X_1,\dots,X_n$ are independent.
\end{proposition*}

\begin{pf}
Without loss of generality, we may, and will, assume that the marginal
distribution of each single $X_i$ is that of Example \ref{arch_ex}
above, that is,\ $P(X_i\leq x)=\psi(-x)$ for $x\leq0$, since we can
always transform each coordinate monotonically so that it has the
proposed distribution after transformation. Thus, with $x_1,x_2,x_3\geq0$,
\begin{eqnarray*}
\Prob(X_3\leq-x_3|X_2\leq-x_2,X_1\leq-x_1)&=&\frac{\Prob(X_3\leq
-x_3,X_2\leq-x_2,X_1\leq-x_1)}{\Prob(X_2\leq-x_2,X_1\leq-x_1)}\\
&=&\frac{\psi(x_1+x_2+x_3)}{\psi(x_1+x_2)},\\
\Prob(X_3\leq-x_3|X_2\leq-x_2)&=&\frac{\Prob(X_3\leq-x_3,X_2\leq
-x_2)}{\Prob(X_2\leq-x_2)}\\
&=&\frac{\psi(x_2+x_3)}{\psi(x_2)},\\
\end{eqnarray*}
so, by the Markov property,
\[
\frac{\psi(x_2+x_3)}{\psi(x_2)}=\frac{\psi(x_1+x_2+x_3)}{\psi(x_1+x_2)}.
\]
Let $f(x)=\psi(x_2+x)/\psi(x_2)$ so that the above equation is
equivalent to
\[
f(x_3) = \frac{f(x_1+x_3)}{f(x_1)}
\]
and thus $f(x_1+x_3)=f(x_1)f(x_3)$, which implies that $f(x) = \mathrm{e}^{-cx}$
for some constant $c$. Putting $x=-x_2$ yields
\[
\mathrm{e}^{cx_2}=f(-x_2)=\frac{\psi(0)}{\psi(x_2)} = \frac{1}{\psi(x_2)}
\]
and thus $\psi(t) = \mathrm{e}^{-ct}$. Hence, $\phi(s)=-\frac{1}{c}\log s$ so
that the copula
\[
\psi\bigl(\phi(u_1)+\cdots+\phi(u_n)\bigr) = u_1\cdots u_n,
\]
that is,\ all $X_1,\dots,X_n$ are independent.
\end{pf}

\section{Idempotent copulas}\label{idem_sec}

A copula is said to be \textit{idempotent} if $C\ast C=C$. In this
section, we will investigate Markov chains with idempotent copulas
whose probabilistic structure will turn out to be quite peculiar.

\begin{example}\label{ord_sum}
Let $I_i=[a_i,b_i]$, $i=1,2,\dots,$ be a set of disjoint intervals in
$[0,1]$. Let $I_0=[0,1]\setminus\bigcup_{i\geq1}I_i$ and let
$p_i=\lambda
(I_i)$ be the Lebesgue measure of each set $I_i$, $i=0,1,\ldots.$
Consider the random variable $(U,V)$ that has the following
distribution: $(U,V)$ is uniform on $I_i\times I_i$ with probability
$p_i$ for $i=1,2,\ldots$ and $U=V$ with $U$ uniform on $I_0$ with
probability $p_0$. Let $D$ be the copula of $(U,V)$. ($D$ is a
so-called ``ordinal sum'' of copies of $\Pi$ and $M$; see
\cite{nelsen}, Chapter~3.2.2.) We have
\[
f_D(x,u) = x\mathbf{1}(x\in I_0)+\sum_{i\geq
1}\bigl((b_i-a_i)u+a_i\bigr)\mathbf{1}(x\in I_i).
\]
It is easy to check that $f_D(f_D(x,u),v)=f_D(x,v)$ so that $D\ast
D=D$, that is,\ $D$ is idempotent. If $U_0,U_1,\dots$ is a Markov chain
governed by the copula $D$, then all $U_0,U_1,\dots$ lie in the same
set $I_{\iota}$, where the random index $\iota$ differs from
realization to realization.
\end{example}

If $C$ is idempotent and $L$ and $R=L^{\mathrm{T}}$ are two copulas satisfying
$L\ast R=M$, then $R\ast C\ast L$ is also idempotent. Darsow \textit{et
al.}\ \cite{dno} conjectured that all idempotent copulas could be
factored in this form with $C$ as in Example \ref{ord_sum}. We will
show that the class of idempotent copulas, even though they correspond
to quite a restricted kind of dependence, is richer than what can be
covered by that characterization. If a Markov chain $X_0,\dots,$ which,
without loss of generality, we assume is in $[0,1]$ is governed by the
copula $R\ast D\ast L$, then all $f_R(X_0,u_0),f_R(X_1,u_1),\dots$ are
in the same set $I_{\iota}$ for some random $\iota$ and there are only
countably many such possible sets. Note that $f_R(x,u)$ is a function
of $x$ only.

We start with some background on spreadable and exchangeable sequences,
with notation from \cite{kallenberg}, which will be useful.

\begin{definition}
An infinite sequence $\xi_1,\xi_2,\dots$ is said to be \emph
{exchangeable} if
\[
(\xi_1,\xi_2,\dots)\stackrel{d}{=}(\xi_{k_1},\xi_{k_2},\dots)
\]
for all permutations $(1,2,\dots)\mapsto(k_1,k_2,\dots)$ which affect
a finite set of numbers. The sequence is said to be \textit{spreadable}
if the equality in distribution is required only for strictly
increasing sequences $k_1<k_2<\cdots.$
\end{definition}

Assuming that the sequence $\xi_1,\xi_2,\dots$ takes its values in a
Borel space, exchangeability and spreadability are, in fact, equivalent
and these notions are also equivalent to the values of the sequence
being conditionally i.i.d.\ with a random distribution $\eta$.
Furthermore, $\eta$ can be recovered from the sequence since $\eta
=\lim
_{n\to\infty}\frac1n\sum_{k\leq n}\delta_{\xi_k}$ a.s., that is,\
$\eta
$ is the almost sure limit of the empirical distribution of $(\xi
_1,\dots,\xi_n)$; see \cite{kallenberg}, Theorem 11.10.

We will need the following observation. If a sequence is conditionally
i.i.d.~given some $\sigma$-algebra $\mathcal{F}$, then $\sigma(\eta
)\subseteq\mathcal{F}$, with $\eta$ as in the previous paragraph.

Let $X_0,X_1,\dots$ be a Markov chain whose Markov copula $C$ is
idempotent. Thus, $C^{\ast n}=C$ for all $n$ and, by the Markov
property, this implies that the sequence is spreadable and\vspace*{-2pt} hence
exchangeable and conditionally i.i.d. Since spreadability implies
$(X_0,X_1)\stackrel{d}{=}(X_0,X_2)$, it is, in fact, equivalent to the
copula being idempotent. This was noted by Darsow \textit{et
al.} \cite{dno}, but we can take the analysis further by using the fact that the
sequence is, in particular, conditionally i.i.d.\ given $\sigma(\eta)$,
where $\eta$ is as above. Thus, for all $n$,
\[
\Prob(X_{n+1}\in\cdot|X_0,\dots,X_n) = \Prob(X_{n+1}\in\cdot
|X_n) =
\Prob(X_{n+1}\in\cdot|X_0),
\]
where the first equality is due to the Markov property and the second
is due to the exchangeability. Therefore,
\begin{eqnarray*}
\Prob\Biggl(\bigcap_{i=0}^n\{X_i\in A_i\}\Biggr)&=& \Prob(X_n\in
A_n|X_{n-1}\in A_{n-1},\dots,X_0\in A_0)\\
&&{}\times\Prob(X_{n-1}\in A_{n-1}|X_{n-2}\in A_{n-2},\ldots,X_0\in
A_0)\cdots\Prob(X_1\in A_1)\\
&=&\Prob(X_n\in A_n |X_0\in A_0)\cdots\Prob(X_1\in A_1 |X_0\in
A_0)
\Prob(X_0\in A_0)
\end{eqnarray*}
and thus $X_1,X_2,\dots$ are conditionally i.i.d.\ given $X_0$, that
is,\ $\sigma(\eta)\subseteq\sigma(X_0)$.

\setcounter{example}{7}
\begin{example}[(Continued)]
Note that $\iota$ is a function of $U_0$ since $\iota$ is the index of
the set that $U_0$ lies in: $\Prob(U_0\in I_{\iota})=1$. It is clear
that the random variables $U_0,U_1,\dots$ that constitute the Markov
chain of the example are i.i.d.~given $U_0$ since all other random
variables are either uniformly distributed on $I_{\iota}$ if $\iota
=1,2,\dots$ or identically equal to $U_0$ if $\iota=0$ (and constant
random variables are independent). Here, the random measure
\[
\eta=\delta_{U_0}\mathbf{1}(\iota=0) + \sum_{i\geq1}\frac
{1}{p_i}\lambda|_{I_i}\mathbf{1}(\iota=i),
\]
where $\delta_x$ is the point mass at $x$ and $\lambda|_{I}$ is the
Lebesgue measure restricted to the set $I$. Since $\eta$ is a function
of $U_0$, we have $\sigma(\eta)\subseteq\sigma(U_0)$.
\end{example}

The following example shows how the proposed characterization fails.

\begin{example}
Let $J_x=\{2^{-n}x, n\in\mathbb{Z}\}\cap(0,1)$ for all $x\in(0,1)$. It
is clear that $\{J_x\}_{x\in[1/2,1)}$ is a partition of $(0,1)$.
Let $m(x)=\max J_x$. We can construct a stationary Markov chain by
letting $U_0$ be uniform on $(0,1)$ and
\[
P\bigl(U_{k+1}=2^{-n}m(x)|U_k=x\bigr)=2^{-(n+1)}
\]
for $n=0,1,2,\dots$ and $k=0,1,\dots.$ Let $E$ be the copula of this
Markov chain. As function $f_E$, we can take
\[
f_E(x,u)=\sum_{n=0}^{\infty}2^{-n}m(x)\mathbf{1}\bigl(2^{-(n+1)}\leq u < 2^{-n}\bigr).
\]
Given $U_0=x$, the rest of the values of the Markov chain
$U_1,U_2,\dots
$ are independent on $J_x$, so the process is conditionally i.i.d. and
$E$ is thus idempotent. Here, the random measure
\[
\eta=\sum_{n=0}^{\infty}2^{-(n+1)}\delta_{2^{-n}m(U_0)},
\]
where $U_0$ is uniform on $(0,1)$. Thus, $\sigma(\eta)\subseteq
\sigma
(U_0)$ is apparent. We note that the cardinality of the set of the
disjoint sets $\{J_x\}_{x\in(\frac12,1]}$ that gives the possible\vspace*{-2pt}
ranges of the Markov chain is uncountable and the copula $E$ can
therefore not be of the form $L^{\mathrm{T}}\ast D\ast L$.
\end{example}

\section*{Acknowledgement} I wish to thank an anonymous referee for
suggesting a stronger statement and a more stringent argument for the
proposition in Section \ref{arch_sec}.

\printhistory

\end{document}